\documentclass[10pt]{article}
\usepackage{amssymb}
\usepackage{amsthm}
\usepackage{epsf,epsfig,amsfonts,amsmath}
\usepackage{graphicx}
\usepackage{a4wide}
\newtheorem{corollary*}{Corollary}

\newcommand{\be}{\begin{equation}}
\newcommand{\ee}{\end{equation}}

\newcommand{\weg}[1]{}   
\begin{document}

\title{On global geodesic mappings of $n$-dimensional surfaces of revolution }
\author{I. Hinterleitner}

\maketitle
\begin{abstract}

In this paper we study geodesic mappings of $n$-dimensional surfaces
of revolution. From the general theory of geodesic mappings of
equidistant spaces we specialize to surfaces of revolution and apply
the obtained formulas to the case of rotational ellipsoids. We prove
that such  $n$-dimensional ellipsoids admit non trivial smooth
geodesic deformations onto $n$-dimensional surfaces of  revolution,
which are generally of a different type.

\end{abstract}



{\bf subclass}: 53B20; 53B21; 53B30; 53C22; 53C25; 53C40
%
\def\a{\alpha}
\def\b{\beta}
\def\G{\Gamma}
\def\epsilon{\varepsilon}
\def\s{\sigma}
\def\vn{$V_n$}
\def\vnn{$\bar V_n$}
\def\nad#1#2{\buildrel{#1} \over{#2}\!\!\strut}
\def\pod#1#2{\mathrel{\mathop{#2}\limits_{#1}}\strut}
\def\tran{; transl.~from }

\section{ Introduction }
We consider geodesic mappings of smooth, closed $n$-dimensional
surfaces of revolution. Such mappings may be constructed with the
aid of embedding into an $(n{+}1)$-dimensional auxiliary space. In
\cite{Mat} the embedding space of an ellipsoid is equipped once with
the euclidian metric, once with a different one. Either metric
induces a metric on the ellipsoid, and by a convenient choice of the
metrics of the embedding space, the induced metrics on the
ellipsoids are geodesically equivalent, i. e. have the same
geodesics.

Another method is to leave the euclidian metric of the embedding
space invariant and to deform the ellipsoid. On each of the deformed
(hyper)surfaces the ambient euclidian metric induces a metric. In
\cite{mkv} a one-parameter family of smooth deformations is
introduced, which maps geodesics to geodesics. In other words, they
are geodesic mappings between different (hyper)surfaces, embedded
into the $n+1$-dimensional euclidian space $E_{n+1}$.

On the other hand, all the deformed surfaces being diffeomorphic to
each other, we can forget the embedding and consider them as one and
the same manifold with the global topology of the n-dimensional
sphere $S_n$. From this point of view our approach generates a
one-parameter family of geodesically equivalent metrics on this
manifold. For an ellipsoid the metrics generated in \cite{Mat} are
special cases of this family.

Concretely we take for simplicity a 2-dimensional rotational
ellipsoid in $E_3$, deform it in a certain way according to
\cite{mkv}, and construct metrics on the deformed surfaces, induced
by the euclidian metric of the embedding space. When we pull back
these metrics to the original ellipsoid, we have a one-parameter
family of geodesically equivalent metrics on one surface. At first
instance, this construction is a realization of geodesic mappings
between different surfaces of revolution, further, when the metrics
of the deformed surfaces are pulled back, the foregoing deformation
gives an intuitive impression of the modified, geodesically
equivalent metrics. The generalization to $n$-dimensional surfaces
of revolution is straightforward.

In the first part (section 2) of the present article we give a brief
introduction into the theory of geodesic mappings of equidistant
spaces. Geodesic mappings are diffeomorphisms $V_n\to \bar V_n$
between $n$-dimensional (pseudo-) Riemannian spaces  that  take any
geodesic of $V_n$ into a geodesic of $\bar V_n$. Equidistant spaces
are manifolds which are foliated by $(n{-}1)$-dimensional subspaces.
In sections 3 and 4  we specialise to $n$-dimensional surfaces of
revolution with rotational symmetry around a certain axis.

In the main part (section 5 and 6) we study rotational ellipsoids
and the class of geodesic mappings introduced in \cite{mkv} as
examples. We find that these mappings act nontrivially and map
ellipsoids onto surfaces of revolution of a different kind, the
explicit description of which is found.

\section{Geodesic mappings of equidistant spaces.}
 General problems of the theory of geodesic mappings have been studied by
T.~Levi-Civita~\cite{le} and many others
\cite{am,ei,kag,mkv,mvh,pe,si,so,st}. It is known that a necessary
and sufficient condition for the existence of geodesic mappings
between $n$-dimensional Riemannian manifolds $V_n(M,g)$ and $\bar
V_n(M,\bar g)$, is $(\bar\nabla-\nabla)_X X=2\nabla_X\psi X$,
$\forall X\in TV_n$, where $\psi$ is a function on $M$ and $\nabla$
and $\bar \nabla$ are the Levi-Civita connections of the metric
tensors $g$ and $\bar g$, respectively. In local coordinates this
``coordinate-free" formula introduced above is equivalent to the
Levi-Civita equation \cite{ei,le,mi,mkv,mvh,pe}
$$\bar \G^h_{ij}(x)=\G^h_{ij}(x)+\delta^h_i\psi_j,$$
 where $\psi_i(x)$  are the components of the gradient $\nabla \psi$, $\G^h_{ij}$ and $\bar\G^h_{ij}$ are the Christoffel symbols of $V_n$ and $\bar V_n$.

Geodesic mappings of equidistant spaces were  studied by N.S. Sinyukov
\cite{sieq,si} and, in the following, by many other authors, for example
\cite{hi, hid, him, mi80r, mi80e, mi85, mi89t, mirot, mi, mkv, mvh, misk, pe, sha,ve}.

Here we present a type of geodesic mappings, introduced in
\cite{mi,si}. Assume the equidistant spaces  $V_n$ and $\bar V_n$,
which have metrics of the following shape  \cite{hi,hid}:

\begin{equation}\label{e1}
{\rm d}s^2=a(w)\,{\rm d}w^2+b(w)\,{\rm d}\sigma^2
\end{equation}
and
\begin{equation}\label{e2}
{\rm d}\bar s^2=\frac{p\,a(w)}{(1+qb(w))^2}\ {\rm d}w^2+\frac{p\,b(w)}{1+qb(w)}\ {\rm d}\sigma^2,
\end{equation}
where \ $a$ \ and \ $b$ \ are differentiable functions of the variable $w$,
 \ $p$ \ and \ $q$ are real parameters, and \
${\rm d}\sigma^2$ is a metric of an $(n{-}1)$-dimensional (pseudo-)
Riemannian space $\tilde V_{n-1}$.

These spaces $V_n$ and $\bar V_n$ have common geodesics, i.e. $V_n
\to \bar V_n$ is a geodesic mapping and the function in the
Levi-Civita equation is
\ $\psi=-\frac 12\ln{\left|1+q\cdot b(w)\right|}$.\, \\
If \ $q\,b'(w)\not\equiv 0$,\,  the mapping is non trivial (i.e.~is
not affine). In \cite{hi,hid,mvh} this is shown by comparing the
Levi-Civita connections of both spaces.

\section{$n$-dimensional surfaces of revolution}

It is well known that surfaces of revolution  and their
$n$-dimensional generalizations are equidistant spaces. Geodesic
mappings of such surfaces are studied in the papers by J. Mike\v s
\cite{mirot,mi,mkv}. The results obtained there  guarantee the
existence of nontrivial geodesic mappings of surfaces of revolution
of class $C^1$, including surfaces which are homeomorphic to the
$n$-dimensional sphere.

Let ${\cal S}_n$  denote a surface of revolution of class $C^1$ in the
$(n{+}1)$-dimensional Euclidian space $E_{n+1}$. In~the cartesian
coordinates system $(x^1,\ldots, x^{n+1})$ the surface ${\cal S}_n$ can be
determined by the conditions:
\begin{equation}\label{e3}
 x^i=r(w)\,u^i\quad (i=1,2,\ldots,n), \qquad x^{n+1}=z(w),
\end{equation}
where $u^i$ are the coordinates of the standard unit sphere
$$
S_{n-1}=\left\{ (u^1, u^2,\ldots, u^n)\,\, |\,\, (u^1)^2+ (u^2)^2+\ldots+(u^n)^2=1 \right\},
$$
$w\in [w_1,w_2]$ is a parameter of horizontal type,
$w_1,w_2$ are fixed values (or $-\infty$ and $+\infty$, respectively),
 $r(w),z(w)\in C^1[w_1,w_2]$, $r(w)>0$ for every $w\in(w_1,w_2)$.

If $r(w_1)=r(w_2)=0$ then ${\cal S}_n$ is homeomorphic to the
$n$-dimensional sphere; the points  corresponding to $w_1$ and $w_2$
are called {\it poles}. In the case $r(w_1)=0$ and $r(w_2)\neq 0$
the surface ${\cal S}_n$ is homeomorphic to $E_n$ or an
$n$-dimensional disk. If $r(w_1)=r(w_2)\neq0$ and $z(w_1)=z(w_2)$,
then the surface ${\cal S}_n$ is homeomorphic to the $n$-dimensional
torus.

In these coordinates the metric on ${\cal S}_n$ has the form

\begin{equation}\label{e4}
{\rm d}s^2=\left( r'^2(w)+z'^2(w)\right){\rm d}w^2+r^2(w)\,{\rm d}\sigma^2,
\end{equation}
where ${\rm d}\sigma^2$ is a metric of the sphere $S_{n-1}$.

{\bf Remark.} In the papers \cite{mirot,mi,mkv}, $r'^2+z'^2=1$ is
assumed, so that the  parameter $w$ denotes the arc length along the
``meridians''. For ${\cal S}_n$ to be smooth at the poles,
corresponding to the parameter values $w_*$, the derivatives of $r$
and $z$ must satisfy $|r'(w_*)|= 1$ and $z'(w_*)=0$.

\section{Geodesic mappings of $n$-dimensional surfaces of revolution}
The metric \eqref{e4} of a surface of revolution appears as a
special case of the metric \eqref{e1} of a general equidistant space
with
\begin{equation}\label{e5}
a(w)=r'^2(w)+z'^2(w) \quad \mbox{and} \quad b(w)=r^2(w).
\end{equation}

In the next section we will construct geodesic mappings to spaces
$\bar V_n$, according to section 2, with a metric of the type,
\begin{equation}\label{e6}
{\rm d}\bar s^2=\bar a(w)\,{\rm d}w^2+\bar b(w)\,{\rm d}\sigma^2,
\end{equation}

\begin{equation}\label{e7}
\bar a(w)=p\ \frac{r'^2+z'^2}{(1+q\,r^2)^2}
\quad \mbox{and}\quad
\bar b(w)=p\ \frac{r^2}{1+q\,r^2}\ ,
\end{equation}
where $p,\,q$ are constants and, evidently $\bar a(w),\, \bar b(w)\neq 0$
and $\bar a(w),\,\, \bar b(w)< \infty$,  when
\begin{equation}\label{e8}
p\neq 0;\quad  r'^2(w)+z'^2(w)\neq 0;\quad  1+qr^2(w)\neq 0 \quad   \forall w\in[w_1,w_2].
\end{equation}

When $q=0$, the space $\bar V_n$ is homothetic (with the coefficient
$p\ \neq 0)$ to the surface ${\cal S}_n$ .

From the formulas \eqref{e6} and  \eqref{e7}  follows that the set
of spaces $\bar V_n$, to which ${\cal S}_n$ maps geodesically,
depends on  two parameters $p$ and $q$.

The  parameter $p$ corresponding to homothetic transformations can
be taken positive, $q$  is restricted by the condition
$1+qr^2(w)\neq 0$ for every $w\in[w_1,w_2]$.

In the case
$$
1+qr^2(w)<0
$$
the metric of $\bar V_n$ has Minkowski signature and this space can
not be realized as surface on a Euclidean space. In this case
$r(w)\neq 0$ for all $w$ and ${\cal S}_n$ is homeomorphic to an
$n$-dimensional cylinder, and
$$\hbox{$q\in\left( -\infty, \pod{w\in[w_1,w_2]}{-\rm max}
\ r^{-2}(w)\right).$}$$

On the other hand if $$1+qr^2(w)>0$$ the metric of $\bar V_n$ is
positive definite and $q$ lies in the interval $q\in\left(
\pod{w\in[w_1,w_2]}{-\rm min} \ r^{-2}(w),+\infty, \right).$
\\
In this case the metric \eqref{e6} might but need not be realised on
some surface of revolution. The ellipsoids considered in the
remainder  of this paper are mapped to closed surfaces of
revolution.

\section{A special class of geodesic mappings}
We suppose that the $n$-dimensional surface of revolution ${\cal
S}_n$, which is defined by equations \eqref{e3}, is  mapped on an
$n$-dimensional surface of revolution $\bar {\cal S}_n$, which is
defined  in a similar manner:
\begin{equation}\label{e9}
x^i=\bar r(w)u^i,\,\,(i=1,2,\ldots,x^n),\quad x^{n+1}=\bar z(w),
\end{equation}
where $\bar r$ and $\bar z$ are functions of the same parameter $w$
as before. In the following we restrict our attention to 2
dimensional surfaces  generated by the rotation of a curve around an
axis, denoted by $z$. The generalization to higher-dimensional
hypersurfaces is straightforward. The curve is formulated in
parametric form $r(w)$, $z(w)$, where $r$ is a radial variable and
the parameter $w$ is the arc length. The relation between $z(w)$ and
some suitably given function $r(w)$ is given by
\begin{equation}\label{e10}
z(w)=\int_{w_1}^we(\tau)\sqrt{1-r'^2(\tau)}\,{\rm d}\tau,
\end{equation}
where $e(\tau)$ can assume the values $\pm1$ on different pieces of
the curve. In the following we assume $e=+1$ throughout the entire
curve, so that $z$ grows monotonically with the length $w$ and the
surface is convex. In \cite{mkv} this case is denoted by the terminus
``simple surface of revolution". $w$ runs from $w_1$ to $w_2$,
$w_2-w_1$ is the total length of the curve. We assume a closed
surface, i.e. $r(w_1)=r(w_2)=0$. For the sake of smoothness at the
rotation poles we assume also $\left|\frac{{\rm d}r}{{\rm
d}w}\right|=1$ at $w_1$ and $w_2$, then $ \frac{{\rm d}z}{{\rm
d}w}=0 $ at the poles follows automatically from \eqref{e10}.

In \cite{mkv} the curve $(r(w),z(w))$ is acted upon by the following
one-parameter family of transformations $f$: ${\cal S}_n\to \bar {\cal S}_n$,
labeled by $a$. It is defined by the following action on the
functions $r(w)$ and $z(w)$:
\begin{equation}\label{e11}
\bar r(w)=\frac{r(w)}{\sqrt{1+ar^2(w)}}, \hspace{1cm} \bar
z(w)=\int_{w_1}^w\sqrt{\frac{1+ar^2(\tau)-r'^2(\tau)}{(1+ar^2(\tau))^3}}\,{\rm
d}\tau.
\end{equation}
The parameter $w$ is the same as before, therefore it is {\em not}
the length parameter of the curve $(\bar r,\bar z)$. The conditions
of smoothness $\bar r=0$, $\frac{{\rm d}\bar r}{{\rm d}w}=\pm1$ and
$\frac{{\rm d}\bar z}{{\rm d}w}=0$ at the poles $w=w_1$ and $w=w_2$
are satisfied, provided they are satisfied for $r$.

It is further shown in \cite{mirot} from comparison of the metrics on the
surfaces generated by rotation of the curves $(r,z)$ and $(\bar
r,\bar z)$ that these mappings leave geodesics invariant. The metric of
${\cal S}_n$ of the form \eqref{e4},
\begin{equation}\label{e12}
{\rm d}s^2={\rm d}w^2+r^2(w)\,{\rm d}\sigma^2,
\end{equation}
is mapped to
\begin{equation}\label{e13}
{\rm d}\bar s^2=\left(\bar r'^2(w)+\bar z'^2(w) \right){\rm
d}w^2+\bar r^2(w){\rm d}\sigma^2.
\end{equation}
In terms of $r$ and $z$ this is explicitly
\begin{equation}\label{e14}
{\rm d}\bar s^2=\frac{{\rm d}w^2}{[1+ar^2(w)]^2}+\frac{r^2(w)\,{\rm
d}\sigma^2}{1+ar^2(w)}.
\end{equation}
The relation between the metrics \eqref{e12} and \eqref{e14} is of
the same type as between \eqref{e1} and \eqref{e2}. In the case of a
2-dimensional surface of revolution embedded into $E_3$, $\sigma$ is
simply the angle around the rotation axis.

>From the point of view of the embedding into 3-dimensional euclidian
space, the transformation $f$ is a smooth deformation. On the other
hand, the formulation \eqref{e14} suggests to consider the original
and the transformed surface as one and the same manifold,
coordinatized by $w$ and the rotation angle $\sigma$, and $f$ as a
transformation of metrics. Then from \eqref{e14} it follows that
both the metrics ${\rm d}s$ and ${\rm d}\bar s$ have the same
geodesics.

\section{Application to rotational ellipsoids}

In the foregoing section we have seen a class of nontrivial geodesic
mappings between smooth surfaces of revolution, which are
homeomorphic to a sphere. Now we take as a concrete example a
rotational ellipsoid, embedded into the 3-dimensional Euclidian
space, and investigate its deformation by the considered geodesic
mappings. This is done in a local coordinate patch, covering one
half of the surface. Rather than in terms of the arc length $w$ we
formulate it in terms of the angular variable $\varphi$,
\begin{equation}\label{el5}
r(\phi)=k\sin\varphi, \hspace{1cm} z(\varphi)=1-\cos\varphi.
\end{equation}
The squared element of the arc length is
\begin{equation}\label{e16}
{\rm d}w^2={\rm d}r^2+{\rm d}z^2=(k^2\cos^2\varphi+\sin^2\varphi)\,{\rm
d}\varphi^2.
\end{equation}
We choose $w_1=w(\varphi=0)=0$, so that the origin of $\varphi$ and the
arc length coincide, then $w_2=w(\varphi=\pi)$ is half of the
circumference of the ellipse. The condition $r(w_1)=r(w_2)=0$ is
fulfilled and
\begin{equation}\label{17}
\frac{{\rm d}r}{{\rm d}w}=\frac{{\rm d}r}{{\rm d}\varphi}\frac{{\rm
d}\varphi}{{\rm
d}w}=\frac{k\,\cos\varphi}{\sqrt{k^2\cos^2\varphi+\sin^2\varphi}},
\end{equation}
so also $\frac{{\rm d}r}{{\rm d}w}(w_1)=1$ and $\frac{{\rm d}r}{{\rm
d}w}(w_2)=-1$ are satisfied.

To carry out the transformation \eqref{e11}, we refer to the parameter
$\phi$ rather than to $w$ for calculational convenience. What is
important, is to use the same parameter for $(r,z)$ and $(\bar
r,\bar z)$. In terms of $\varphi$ we have
\begin{equation}\label{e18}
\bar r(\varphi)=\frac{k\,\sin\varphi}{\sqrt{1+ak^2\sin^2\varphi}}
\end{equation}
and
\begin{equation}\label{e19}
\bar
z(\varphi)=\int_0^\varphi\sqrt{\frac{1+ar^2(\varphi')-r'^2(\varphi')}{(1+ar^2(\varphi'))^3}}\,\frac{{\rm
d}w}{{\rm d}\varphi'}\,{\rm d}\varphi'.
\end{equation}
Note that here and in the following $r'$ means always the derivative
with respect to $w$, even when written as function of $w$, so
$r'(\varphi)$ is $\frac{{\rm d}r}{{\rm d}\varphi}\,\frac{{\rm d}\varphi}{{\rm
d}w}$ as function of $\varphi$.

Explicitly we find
\begin{equation}\label{e20}
\bar
r'(\varphi)=\frac{k\,\cos\varphi}{(k^2\cos^2\varphi+\sin^2\varphi)^{\frac{1}{2}}
(1+ak^2\sin^2\varphi)^\frac{3}{2}},
\end{equation}
the maximal value of $\bar r$, $\bar r_{\rm
max}=\frac{k}{\sqrt{1+ak^2}}$, occurs at $\varphi=\frac{\pi}{2}$, like
for the original ellipsoid.

To find out what kind of curve is $(\bar r(\varphi),\bar z(\varphi))$, we
should solve the integral \eqref{e20} explicitly, which is
complicated. So we consider rather the derivative $\frac{{\rm d}\bar
z}{{\rm d}\bar r}$, which gives a differential equation for the
curve, and eliminate the parameter. This is done in several steps:

First we express $\frac{{\rm d}\bar z}{{\rm d}\bar r}$ in the form
$\frac{{\rm d}\bar z}{{\rm d}\varphi}/\frac{{\rm d}\bar r}{{\rm
d}\varphi}$. From \eqref{e18} and \eqref{e19} we get
$$\frac{{\rm d}\bar r}{{\rm
d}\varphi}=\frac{k\cos\varphi}{(1+ak^2\sin^2\varphi)^{\frac{3}{2}}}$$ and
$$\frac{{\rm d}\bar z}{{\rm d}\varphi}=\sqrt{\frac{1+ar^2(\varphi)-r'^2(\varphi)}{(1+ar^2(\varphi))^3}}\,
\frac{{\rm d}w}{{\rm d}\varphi}.$$

Then from the definitions \eqref{el5} and the explicit equation of
the ellipse
$$(1-z)^2+\frac{r^2}{k^2}=1$$
we express $\sin\varphi$ and $\cos\varphi$ in terms of $r$ and find
\begin{equation}\label{e21}
\frac{{\rm d}\bar z}{{\rm d}\bar
r}=\frac{r\sqrt{1+ak^4+ar^2-a^2k^2r^2}}{k\sqrt{k^2-r^2}}
\end{equation}
in terms of $r$.

Now we insert the inverse of \eqref{e18},
\begin{equation}
r=\frac{\bar r}{\sqrt{1-a\bar r^2}}
\end{equation}
to express this derivative in terms of $\bar r$,
\begin{equation}
\frac{{\rm d}\bar z}{{\rm d}\bar r}=\frac{\bar
r\,\sqrt{\frac{1}{k^2}+ak^2-a(1+ak^2)\bar r^2}}{\sqrt{1-a\bar
r^2}\sqrt{k^2-(1+ak^2)\bar r^2}}.
\end{equation}
At last, for a direct comparison with the corresponding differential
equation for an ellipse,
\begin{equation}
\frac{{\rm d}z}{{\rm d}r}=\frac{r}{k\sqrt{k^2-r^2}},
\end{equation}
we carry out a scale transformation
\begin{equation}
\hat r=\bar r\sqrt{1+ak^2}, \hspace{1cm}\hat z=\bar
z\sqrt{1+ak^2},
\end{equation}
so that the maximal value of $\hat r$ is equal to $k$, like the maximal
value of $r$ in the case of the ellipse and the radial extensions of
both surfaces are the same. In terms of these variables, finally,
\begin{equation}\label{barbar}
\frac{{\rm d}\hat z}{{\rm d}\hat r}=\frac{\hat
r}{k\sqrt{k^2-\hat r^2}}\cdot
\sqrt{\frac{1+ak^2(k^2-\hat r^2)}{1+a(k^2-\hat r^2)}}.
\end{equation}
>From this we can see that the transformed curve is of a different
type than an ellipse, the differential equation of the corresponding
ellipse is modified by the factor right to the dot in
\eqref{barbar}. At the maximal values of the radial variables, i. e.
at the ``equator'', both the derivatives $\frac{{\rm d}z}{{\rm d}r}$
for the ellipse and $\frac{{\rm d}\hat z}{{\rm d}\hat r}$ for
the deformed curve go to infinity, corresponding to the fact that
$r$ and $\sigma$ provide only a local chart for one half of the
surface.

An interesting feature of these transformations is that they leave
circles ($k=1$) invariant (up to a scale factor $\sqrt{1+a}$), in
the limit of a large transformation parameter $a$ the modification
factor in \eqref{barbar} goes to $k$ and the (rescaled) transformed curve approaches is
a circle.

The metric of the resulting surface of revolution,
$$
{\rm d}s^2=\left(1+\frac{{\rm d}\hat z^2}{{\rm d}\hat
r^2}\right){\rm d}\hat r^2+\hat r^2{\rm d}\sigma^2,
$$
is
\begin{equation}\label{metric}
{\rm d}s^2=\frac{k^2+ak^4+\left(\frac{1}{k^2}-ak^2-1\right)\hat
r^2}{(k^2-\hat r^2)(1+ak^2-a\hat r^2)}\,{\rm d}\hat r^2+\hat
r^2\,{\rm d}\sigma^2.
\end{equation}
This form of the metric in terms of $\hat r$ is local and applies
only to the lower or the upper half of the surface.
It can be generalized without problems to higher dimensions, when
the circles with constant $\hat{z}$ are replaced by
higher-dimensional spheres. Then ${\rm d}\sigma$ has only to be
replaced by the solid angle element ${\rm d}\Omega$ of the
corresponding dimension.

This metric can be pulled back to the original ellipsoid by simply
expressing $\hat r$ in terms of $r$,
\begin{equation}\label{e28}
{\rm d}s^2=(1+a
k^2)\left[\frac{k^2+\left(\frac{1}{k^2}-1\right)r^2}{(k^2-r^2)(1+a
r^2)^2}\,{\rm d}r^2+\frac{r^2}{1+a r^2}\,{\rm
d}\sigma^2\right],
\end{equation}
whereas the metric on the original ellipsoid is
\begin{equation}\label{e29}
{\rm d}s^2=\frac{k^2+(\frac{1}{k^2}-1)r^2}{k^2-r^2}\,{\rm
d}r^2+r^2\,{\rm d}\sigma^2.
\end{equation}

For an explicit expression of the deformed surfaces, we calculate
the equations of the ``meridians" in the form $\hat z(r)$, where
$\hat z$ and $\hat r$ are cartesian coordinates of a cross-section
through the rotation axis. For this purpose we integrate
\eqref{barbar}, from now on we drop the hats on $r$ and $z$.
 We begin with the substitution
\begin{equation}
\sin^2\phi=\frac{k^2-r^2}{k^2-\frac{1}{a}-r^2}.
\end{equation}
Then
\begin{equation}\label{e31}
z(r)=-\frac{1}{k\sqrt{a}}\int_{\phi(0)}^{\phi(r)}\frac{\sqrt{1-(1-k^2)\sin^2\phi}}
{\cos^2\phi}\,{\rm d}\phi,
\end{equation}
where
$$\phi(0)=\arcsin\sqrt{\frac{ak^2}{1+ak^2}}\hspace{5mm}\mbox{and}\hspace{5mm}
\phi(r)=\arcsin\sqrt{\frac{a(k^2-r^2)}{1+ak^2-ar^2}}.$$
 Integrating \eqref{e31} by parts gives
$$\begin{array}{l}
\displaystyle\int_{\phi(0)}^{\phi(r)}\sqrt{1-(1-k^2)\sin^2\!\phi}\:\frac{{\rm
d}\phi}{\cos^2\phi}=\left.\sqrt{1-(1-k^2)\sin^2\!\phi}\:\tan\phi\right|_{\phi(0)}^{\phi(r)}\\[3mm]
\hspace{5mm}-\displaystyle\frac{1}{(1-k^2)}\int_{\phi(0)}^{\phi(r)}\sqrt{1-(1-k^2)\sin^2\!\phi}\:{\rm
d}\phi+\frac{1}{(1-k^2)}\int_{\phi(0)}^{\phi(r)}\frac{{\rm
d}\phi}{\sqrt{1-(1-k^2)\sin^2\!\phi}}\end{array}$$ where the last
two integrals are the standard elliptic integrals of the second and
first kind \cite{dw}
$$E(\phi,\kappa)=\int_0^\phi\sqrt{1-\kappa^2\sin^2\!\phi}\;{\rm
d}\phi \hspace{3mm}\mbox{and}\hspace{3mm}
F(\phi,\kappa)=\int_0^\phi\frac{{\rm
d}\phi}{\sqrt{1-\kappa^2\sin^2\!\phi}}.$$ Inserting back $r$ gives
finally
\begin{equation}\label{e32}
\begin{array}{l}
z(r)=-\displaystyle\frac{\sqrt{k^2-r^2}}{k}\sqrt{\frac{1+ak^4-ak^2r^2}{1+ak^2-ar^2}}+
\sqrt{\frac{1+ak^4}{1+ak^2}}+\\[6mm]
\hspace{4mm}\displaystyle\frac{1}{\sqrt{a}\,k(1-k^2)}\left[E\left(\arcsin\sqrt{\frac{k^2-r^2}
{k^2+\frac{1}{a}-r^2}},\sqrt{1-k^2}\right)-E\left(\arcsin\frac{k}{\sqrt{k^2+\frac{1}{a}}},
\sqrt{1-k^2}\right) \right.\\[5mm]
\left.\hspace{4mm}-F\left(\arcsin\displaystyle\sqrt{\frac{k^2-r^2}{k^2+\frac{1}{a}-r^2}},\sqrt{1-k^2}\right)
+F\left(\arcsin\displaystyle\frac{k}{\sqrt{k^2+\frac{1}{a}}},\sqrt{1-k^2}\right)\right]
\end{array}
\end{equation}

For small values of the parameter $a$ we have in linear approximation
$$
z(r)=1-\frac{\sqrt{k^2-r^2}}{k}-a\cdot\frac {k^2(1-k^2)}{6}+a\cdot \frac {1-k^2}{6k}(k^2-r^2)^{\frac 32}+O(a^2).
$$

From this we see that the $z$ coordinate of the ``equator", where
$r=k$, is shifted by $-a\cdot \frac {k^2(1-k^2)}{6}$, so according
to the sign of $a$ and $1-k^2$ the surface becomes ``compressed" or
``elongated".

\section{Summary}

We have considered two aspects of geodesic mappings of ellipsoids.
\eqref{barbar} and \eqref{e32} describe the transformations as
deformations in $E_3$. An interesting property is that on a sphere
as a special case of an ellipsoid these transformations act as
identity, whereas they act highly non trivially on general
ellipsoids. In the limit of large transformation parameters the
transformed surfaces approach a sphere as limiting surface.

The second aspect, represented by \eqref{e28} and \eqref{e29},
concerns geodesic transformations of the metric on a manifold
homeomorphic to the sphere, in a accordance with \cite{voss}, where
it is shown by application of a classical theorem by Dini \cite{di}
that there is (up to homothety) a one-parameter family of
geodesically equivalent metrics on $S_2$.

{\bf Acknowledgments.} This work was partially supported by the
Ministry of Education, Youth and Sports of the Czech Republic, Research \&\
Development, project No. 0021630511.\\

Irena Hinterleitner\\
Brno University of Technology, Faculty of Civil Engineering, Dept.
of Math., \v Zi\v zkova 17,\\
 602 00 Brno, Czech Rep.\\
{\it E-mail}: hinterleitner.irena@seznam.cz

\end{document}